\documentclass[12pt]{article}

\newcommand{\qed}{{\hfill\rule{4pt}{7pt}}\medskip}

\usepackage{amsmath, epsfig}
\usepackage{amssymb}
\usepackage{amsfonts}
\usepackage{latexsym}

\newtheorem{thm}{Theorem}

\newtheorem{cor}[thm]{Corollary}

\newtheorem{lem}[thm]{Lemma}

\def\lcm{{\rm lcm}}

\def\pf{\noindent {\it Proof.} }

\begin{document}

\begin{center}
{\Large\bf On the least common multiple of $q$-binomial coefficients}
\end{center}

\vskip 2.5mm
\centerline{\large Victor J. W. Guo}

\vskip 2mm
\centerline{\footnotesize
Department of Mathematics, East China
Normal University, Shanghai 200062, People's Republic of China}

\begin{center}
{\footnotesize\tt jwguo@math.ecnu.edu.cn
}
\end{center}

\vskip 0.7cm \noindent{\bf Abstract.} In this paper, we prove the following identity
\begin{align*}
\lcm\left({n\brack 0}_q,{n\brack 1}_q,\ldots,{n\brack n}_q\right)
=\frac{\lcm([1]_q,[2]_q,\ldots,[n+1]_q)}{[n+1]_q},
\end{align*}
where ${n\brack k}_q$ denotes the $q$-binomial coefficient and $[n]_q=\frac{1-q^n}{1-q}$.
This result is a $q$-analogue of an identity of Farhi [Amer. Math. Monthly, November (2009)].

\vskip 0.2cm
\noindent{\it Keywords}: least common multiple; $q$-binomial coefficient; cyclotomic polynomial

\vskip 0.5cm
\noindent{\it AMS Subject Classifications (2000)}: 11A07; 05A30


\section{Introduction}
An equivalent form of the prime number theorem states that
$\log \lcm(1,2,\ldots,n)\sim n$ as $n\to\infty$ (see, for example, \cite{HW}).
Nair \cite{Nair} gave a nice proof for the well-known estimate $\lcm\{1,2,\ldots,n\}\geq 2^{n-1}$,
while Hanson \cite{Hanson} already obtained $\lcm\{1,2,\ldots,n\}\leq 3^{n}$.
Recently, Farhi \cite{Farhi} established the following interesting result.
\begin{thm}[Farhi]\label{thm:farhi}
For any positive integer $n$, there holds
\begin{align}
\lcm\left({n\choose 0},{n\choose 1},\ldots,{n\choose n}\right)
=\frac{\lcm(1,2,\ldots,n+1)}{n+1}.  \label{eq:farhi}
\end{align}
\end{thm}
As an application, Farhi shows that $\lcm\{1,2,\ldots,n\}\geq 2^{n-1}$
follows immediately from \eqref{eq:farhi}.

The purpose of this note is to give a $q$-analogue of \eqref{eq:farhi} by using cyclotomic polynomials.
Recall that a natural $q$-analogue of the nonnegative integer $n$ is given by $[n]_q=\frac{1-q^n}{1-q}$.
The corresponding $q$-factorial is $[n]_q!=\prod_{k=1}^n[k]_q$ and the $q$-binomial coefficient
${M\brack N}_q$ is defined as
\[
{M\brack N}_q=
\begin{cases}
\displaystyle\frac{[M]_q!}{[N]_q! [M-N]_q!}, &\text{if $0\leq N\leq M$},\\[5pt]
0,&\text{otherwise.}
\end{cases}
\]

Let $\lcm$ also denote the least common multiple of a sequence of polynomials in $\mathbb{Z}[q]$.
Our main result can be stated as follows:
\begin{thm}\label{thm:main}
For any positive integer $n$, there holds
\begin{align}
\lcm\left({n\brack 0}_q,{n\brack 1}_q,\ldots,{n\brack n}_q\right)
=\frac{\lcm([1]_q,[2]_q,\ldots,[n+1]_q)}{[n+1]_q}.  \label{eq:main}
\end{align}
\end{thm}

\section{Proof of Theorem \ref{thm:main}}

Let $\Phi_n(x)$ be the $n$-th {\it cyclotomic polynomial}.
The following easily proved result can be found in \cite[(10)]{KW} and
\cite{GZ}.

\begin{lem}\label{prop:factor}
The q-binomial coefficient ${n\brack k}_q$ can be factorized into
$$
{n\brack k}_q=\prod_{d}\Phi_d(q),
$$
where the product is over all positive integers $d\leq n$ such that
$\lfloor k/d\rfloor+\lfloor (n-k)/d\rfloor<\lfloor n/d\rfloor$.
\end{lem}

\begin{lem}\label{prop:solve-d}
Let $n$ and $d$ be two positive integers with $n\geq d$. Then there exists at least one positive integer $k$
such that
\begin{align}
\lfloor k/d\rfloor+\lfloor (n-k)/d\rfloor<\lfloor n/d\rfloor  \label{eq:kdfloor}
\end{align}
if and only if $d$ does not divide $n+1$.
\end{lem}
\pf Suppose that \eqref{eq:kdfloor} holds for some positive integer $k$. Let
$$
k\equiv a \pmod d,\qquad (n-k) \equiv b \pmod d
$$
for some $1\leq a,b\leq d-1$. Then $n\equiv a+b\pmod d$ and $d\leq a+b\leq 2d-2$. Namely,
$n+1\equiv a+b+1\not\equiv 0\pmod d$. Conversely, suppose that $n+1\equiv c\pmod d$ for
some $1\leq c\leq d-1$. Then $k=c$ satisfies \eqref{eq:kdfloor}.
This completes the proof. \qed

\noindent{\it Proof of Theorem \ref{thm:main}. } By Lemma \ref{prop:factor}, we have
\begin{align}
\lcm\left({n\brack 0}_q,{n\brack 1}_q,\ldots,{n\brack n}_q\right)
=\prod_{d}\Phi_d(q), \label{eq:d}
\end{align}
where the product is over all positive integers $d\leq n$ such that for some $k$ ($1\leq k\leq n$)
there holds $\lfloor k/d\rfloor+\lfloor (n-k)/d\rfloor<\lfloor n/d\rfloor$.
On the other hand, since
$$
[k]_q=\frac{q^k-1}{q-1}=\prod_{d\mid k,\ d>1}\Phi_d(q),
$$
we have
\begin{align}
\frac{\lcm([1]_q,[2]_q,\ldots,[n+1]_q)}{[n+1]_q}=\prod_{d\leq n,\ d\nmid (n+1)}\Phi_d(q).
\label{eq:d2}
\end{align}
By Lemma \ref{prop:solve-d}, one sees that the right-hand sides of \eqref{eq:d} and \eqref{eq:d2}
are equal. This proves the theorem.  \qed

\section{Theorem \ref{thm:main} is a $q$-analogue of Theorem \ref{thm:farhi} }
In this section we will show that
\begin{align}
\lim_{q\to 1}\lcm\left({n\brack 0}_q,{n\brack 1}_q,\ldots,{n\brack n}_q\right)
=\lcm\left({n\choose 0},{n\choose 1},\ldots,{n\choose n}\right),  \label{eq:limitq=1}
\end{align}
and
\begin{align}
\lim_{q\to 1}\frac{\lcm([1]_q,[2]_q,\ldots,[n+1]_q)}{[n+1]_q}
=\frac{\lcm(1,2,\ldots,n+1)}{n+1}.  \label{eq:limit2}
\end{align}
We need the following property.
\begin{lem}\label{lem:q=1}
For any positive integer $n$, there holds
\begin{align*}
\Phi_n(1)=
\begin{cases}
p,&\text{if $n=p^r$ is a prime power},\\
1, &\text{otherwise.}
\end{cases}
\end{align*}
\end{lem}
\pf See for example \cite[p.~160]{Nagell}.
\qed

In view of \eqref{eq:d}, we have
\begin{align}
\lim_{q\to 1}\lcm\left({n\brack 0}_q,{n\brack 1}_q,\ldots,{n\brack n}_q\right)
=\prod_{d}\Phi_d(1), \label{eq:dq=1}
\end{align}
where the product is over all positive integers $d\leq n$ such that for some $k$ ($1\leq k\leq n$)
there holds $\lfloor k/d\rfloor+\lfloor (n-k)/d\rfloor<\lfloor n/d\rfloor$. By Lemma \ref{lem:q=1},
the right-hand side of \eqref{eq:dq=1} can be written as
\begin{align}
\prod_{\text{primes $p\leq n$}} p^{\sum_{r=1}^{\infty}\max_{0\leq k\leq n}
\left\{\lfloor n/p^r\rfloor-\lfloor k/p^r\rfloor-\lfloor (n-k)/p^r\rfloor\right\}}.  \label{eq:primep}
\end{align}

We now claim that
\begin{align}
&\hskip -3mm
\sum_{r=1}^{\infty}\max_{0\leq k\leq n}
\left\{\lfloor n/p^r\rfloor-\lfloor k/p^r\rfloor-\lfloor (n-k)/p^r\rfloor\right\} \nonumber\\
&=\max_{0\leq k\leq n}\sum_{r=1}^{\infty}
\left(\lfloor n/p^r\rfloor-\lfloor k/p^r\rfloor-\lfloor (n-k)/p^r\rfloor\right).  \label{eq:max}
\end{align}
Let $n=\sum_{i=0}^M a_i p^i$, where $0\leq a_0,a_1,\ldots,a_M\leq p-1$ and
$a_M\neq 0$. By Lemma \ref{prop:solve-d}, the left-hand side of \eqref{eq:max} (denoted LHS\eqref{eq:max})
is equal to
the number of $r$'s such that $p^r\leq n$ and $p^r\nmid n+1$. It follows that
$$
LHS\eqref{eq:max}=
\begin{cases}
0,&\text{if $n=p^{M+1}-1$},\\
M-\min\{i\colon a_i\neq p-1\}, &\text{otherwise.}
\end{cases}
$$
It is clear that the right-hand side of \eqref{eq:max} (denoted RHS\eqref{eq:max}) is less than
or equal to LHS\eqref{eq:max}. If $n=p^{M+1}-1$, then both sides of \eqref{eq:max} are equal to $0$.
Assume that $n\neq p^{M+1}-1$ and
$i_0=\min\{i\colon a_i\neq p-1\}$.
Taking $k=p^{M}-1$, we have
$$
\lfloor n/p^r\rfloor-\lfloor k/p^r\rfloor-\lfloor (n-k)/p^r\rfloor
=\begin{cases}
0,&\text{if $r=1,\ldots,i_0$},\\
1,&\text{if $r=i_0+1,\ldots,M$,}
\end{cases}
$$
and so
$$
\sum_{r=1}^\infty \lfloor n/p^r\rfloor-\lfloor k/p^r\rfloor-\lfloor (n-k)/p^r\rfloor=M-i_0.
$$
Thus \eqref{eq:max} holds. Namely, the expression \eqref{eq:primep} is equal to
\begin{align*}
\prod_{\text{primes $p\leq n$}} p^{\max_{0\leq k\leq n}\sum_{r=1}^{\infty}
\left(\lfloor n/p^r\rfloor-\lfloor k/p^r\rfloor-\lfloor (n-k)/p^r\rfloor\right)}
=\lcm\left({n\choose 0},{n\choose 1},\ldots,{n\choose n}\right).
\end{align*}
This proves \eqref{eq:limitq=1}. To prove \eqref{eq:limit2}, we apply \eqref{eq:d2} to get
\begin{align*}
\lim_{q\to 1}\frac{\lcm([1]_q,[2]_q,\ldots,[n+1]_q)}{[n+1]_q}=\prod_{d\leq n,\ d\nmid (n+1)}\Phi_d(1),
\end{align*}
which, by Lemma \ref{lem:q=1}, is clearly equal to
$$
\frac{\lcm(1,2,\ldots,n+1)}{n+1}.
$$

Finally, we mention that \eqref{eq:max} has the following interesting conclusion.
\begin{cor}
Let $p$ be a prime number and let $k_1,k_2,\ldots, k_m\leq n$, $r_1<r_2<\cdots<r_m$ be positive integers such that
\begin{align*}
\lfloor n/p^{r_i}\rfloor-\lfloor k_i/p^{r_i}\rfloor-\lfloor (n-k_i)/p^{r_i}\rfloor=1 \quad for\ i=1,2,\ldots,m.
\end{align*}
Then there exists a positive integer $k\leq n$ such that
\begin{align*}
\lfloor n/p^{r_i}\rfloor-\lfloor k/p^{r_i}\rfloor-\lfloor (n-k)/p^{r_i}\rfloor=1 \quad for\ i=1,2,\ldots,m.
\end{align*}
\end{cor}


\renewcommand{\baselinestretch}{1}

\end{document}